\documentclass[12pt]{article}

\usepackage{amssymb}
\usepackage{amsmath}
\usepackage{latexsym}
\usepackage{mathrsfs}

\usepackage{color,fancybox}
\newtheorem{theorem}{Theorem}

\newtheorem{definition}[theorem]{Definition}

\newcommand{\po}{{\partial\Omega}}
\newcommand{\nl}{\newline}

\newcommand{\dom}{{\rm Dom}}

\newcommand{\N}{{\mathbb{N}}}

\newcommand{\R}{{\mathbb{R}}}

\newcommand{\cC}{{\cal C}}

\newcommand{\eps}{\epsilon}

\newcommand{\parder}[2]{\frac{\partial{#1}}{\partial{#2}}}

\newcommand{\diver}{{\rm div}}

\newcommand{\darr}[4]{{\left\{\begin{array}{ll}
  {#1}&{#2}\\[0.3cm]
  {#3}&{#4}
\end{array}\right.}}

\newcommand\triplenorm[1]{\left|\hspace{-1pt}\left|\hspace{-1pt}\left|{#1}\right|\hspace{-1pt}\right|\hspace{-1pt}\right| }

\def\Frac{\displaystyle\frac}

\title{Stability estimates in $H^1_0$ for solutions of elliptic equations in varying domains}
\author{Jos\'{e} M. Arrieta \ \  \ Gerassimos Barbatis}
\date{}
\begin{document}
\maketitle

\begin{abstract}
\noindent
We consider second-order uniformly elliptic operators subject to Dirichlet boundary conditions. Such operators are considered on a bounded domain $\Omega$ and on the domain $\phi(\Omega)$ resulting from $\Omega$ by means of a bi-Lipschitz map $\phi$. We consider the solutions $u$ and $\tilde u$ of the corresponding elliptic equations with the same right-hand side $f\in L^2(\Omega\cup\phi(\Omega))$. Under certain assumptions we estimate the difference $\|\nabla\tilde u-\nabla u\|_{L^2(\Omega\cup\phi(\Omega))}$ in terms of certain measure of vicinity of $\phi$ to the identity map. For domains within a certain class this provides estimates in terms of the Lebesgue measure of the symmetric difference of $\phi(\Omega)$ and $\Omega$, that is $|\phi(\Omega)\triangle \Omega|$.
We provide an example which shows that the estimates obtained are in a certain sense sharp.
\end{abstract}

\vspace{11pt}


\section{Introduction}

Let $\Omega$ be a bounded domain in $\R^N$ and let $A=\{A_{ij}(x)\}$ be a measurable positive-definite symmetric
 matrix-valued function on $\R^N$ bounded away from zero and infinity. Let
\[
Lu =- \sum_{i,j=1}^N\parder{}{x_i}\Big\{ A_{ij}(x)\parder{u}{x_j}\Big\}, \qquad \mbox{ on }\Omega,
\]
subject to Dirichlet boundary conditions. Let $\phi:\Omega\to\phi(\Omega)$ be a bi-Lipschitz map and let $\tilde L$ be the analogous operator on
$L^2(\phi(\Omega))$. We fix $f\in L^2(\Omega\cup\phi(\Omega))$ and consider the functions $u\in H^1_0(\Omega)$ and $\tilde u\in H^1_0(\phi(\Omega))$ defined by
\[
Lu =f, \quad \mbox{ in }\Omega \; , \qquad  \tilde L\tilde u =f\, , \quad  \mbox{ in }\phi(\Omega).
\]
Our aim in this article is to estimate $\|\nabla \tilde u-\nabla u\|_{L^2(\Omega\cup\phi(\Omega))}$ in terms of a certain measure of vicinity of $\phi$ to the identity map.

Domain perturbation problems are an important branch of the theory of PDEs. Within spectral theory
one is typically intersted on the stability of the eigenvalues or eigenfunctions of differential operators. There are several recent results on this type of problems; we refer to the article \cite{BL2} and references therein for more on recent progress on domain perturbation problems in spectral theory.

The problem we are interested in involves the stability of the solution $u$ of the equation $Lu=f$ under perturbation of the domain $\Omega$.
In the article \cite{SS} Savar\'{e} and Schimperna obtain very detailed sharp estimates on the variation in $H^1_0$ of $u$ for operators with Lipschitz continuous coefficients. In \cite{BBL} estimates are obtained on the variation of $u$ in $L^2$ for operators with measurable coefficients. See also the articles \cite{AD1,AD2,B} where relevant results were obtained.

In the present article we consider uniformly elliptic operators with measurable coefficients and we prove stability estimates in $H^1_0$ for the solution $u$ of $Lu=f$. A simple example shows that the estimates are in a certain sense sharp. Our main assumptions are, roughly,  that $\Omega$ is perturbed by a global bi-Lipschitz map $\phi$ and that $f$ and $\nabla u$ belong in $L^q$ for some $q>2$. We note that if $f\in L^q$ and $\Omega$ has Lipschitz boundary in the sense that there exists a bi-Lipschitz transformation that maps $\Omega$ onto a domain with $C^1$ boundary, then $\nabla u\in L^q$ by a well known result of Meyers \cite{Me}, so our results are applicable. The proof relies on the so-called pull-back method. The operator $\tilde L$ that acts on $\phi(\Omega)$ induces naturally an operator on $\Omega$, and it is that operator which is then compared with $H$. Hence, in the first section we prove a stability estimate for operators acting on the same space under variation of the coefficients.

The method of proof easily generalizes to other kinds of differential operators, such as higher-order operators, operators subject to Neumann boundary conditions or operators acting on Riemannan manifolds. For the sake of simplicity and brevity we restrict our attention to second-order Dirichlet operators on bounded Euclidean domains.

\section{A general stability theorem}
\label{sec2}

In this section we prove an auxiliary result which we believe is of independent interest.
Let $\Omega$ be a bounded domain in $\R^N$ and for $\epsilon\geq 0$ let $A_{\epsilon}=(A_{\epsilon}^{ij})$ be a family of real, symmetric, matrix-valued measurable functions on $\Omega$ satisfying
\begin{equation}
\frac{1}{c} |\xi|^2 \leq \sum_{i,j}A_{\epsilon}^{ij}(x)\xi_i\xi_j \leq c |\xi|^2 \; , \quad \epsilon\geq 0 , \;\; x\in\Omega , \;\; \xi\in\R^n\, .
\label{0}
\end{equation}
For $\epsilon\geq 0$ we define the self-adjoint operator
\[
H_{\epsilon} =-\sum_{i,j}\parder{}{x_j}\Big( A_{\epsilon}^{ij}(x)\parder{}{x_i} \Big) \; , \qquad \mbox{ on }L^2(\Omega),
\]
subject to Dirichlet boundary conditions on $\partial\Omega$. We now fix $f\in L^2(\Omega)$ and we denote by $u_{\epsilon}$ the solution of $H_{\epsilon}u_{\epsilon}=f$, $\epsilon\geq 0$. Using a standard argument we obtain from (\ref{0}) that $\|u_{\epsilon}\|_{H^1(\Omega)}\leq c$ for all $\epsilon>0$ hence, up to a subsequence, $(u_{\epsilon})$ converges weakly in $H^1_0(\Omega)$.

We now assume that $\|A_\epsilon-A_0\|_{L^p(\Omega)}\equiv \sup_{i,j}\{\|A^{ij}_{\epsilon}-A^{ij}_0\|_{L^p(\Omega)}\} \to 0$ for some (equivalently, for all) $1\leq p<\infty$.
We then easily deduce that the weak $H^1$-limit of $(u_{\epsilon})$ is precisely $u_0$. This implies in particular that the full sequence $(u_{\epsilon})$ converges to $u_0$ weakly in $H^1_0(\Omega)$.
It is a natural question to seek conditions under which the convergence $u_{\epsilon}\to u_0$ is strong in $H^1_0(\Omega)$. To our knowledge this problem has not been studied. In the next theorem we provide four conditions each one of which guarantees strong convergence in $H^1_0(\Omega)$. The first of these conditions also provides an estimate for the rate of convergence $u_{\epsilon}\to u_0$ in $H^1_0(\Omega)$, and is the one that will be used for the domain perturbation problem in the next section.

Given a real symmetric matrix $A$, the matrix $A_+$ is defined by means of the spectral theorem: if $A=\sum \lambda_n (e_n\otimes e_n)$ is the spectral representation of $A$, then $A_+:= \sum (\lambda_n)_+ (e_n\otimes e_n)$.

\begin{theorem}
Assume that any one of the following four conditions is satisfied:
\begin{eqnarray*}
(1) && \mbox{$\nabla u_0\in L^q(\Omega)$ for some $q>2$.}\\[0.12cm]
(2) && \mbox{The eigenfunctions $\{\phi_n\}$ of $H_0$ satisfy $\|\nabla\phi_n\|_{L^q(\Omega)}\leq c\lambda_n^{\gamma}$ for}\\[0.04cm]
&& \mbox{some $q>2$, $\gamma>0$ and all $n\in\N$.}\\[0.12cm]
(3) && \mbox{There exists a compact set $K\subset\Omega$ such that $A_{\epsilon}\to A_0$ in $L^{\infty}(\Omega\setminus K)$.}\\[0.12cm]
(4) && \mbox{$(A_0-A_{\epsilon})_+ \to 0 $ in 
$L^{\infty}(\Omega)$.}\\
\end{eqnarray*}
Then $u_{\epsilon}\to u_0$ in $H^1_0(\Omega)$ as $\epsilon \to 0$. Moreover, in case (1) we have the estimate
\begin{equation}
\|\nabla u_{\epsilon}-\nabla u_0\|_{L^2(\Omega)} \leq c \|\nabla u_0\|_{L^{q}(\Omega)}\|a_{\epsilon}-a_0\|_{L^{2q/(q-2)}(\Omega)} \; , \qquad \epsilon >0.
\label{est}
\end{equation}
\label{thm}
\end{theorem}
{\em Proof of Theorem \ref{thm}.} {\em Part (1).}  Using the standard notation of repeated indices, we have
\begin{eqnarray*}
&&\int_{\Omega}A^{ij}_{\epsilon}\frac{\partial u_{\epsilon}}{\partial x_i}\cdot (\frac{\partial u_{\epsilon}}{\partial x_j}-\frac{\partial u_0}{\partial x_j})dx=\int_{\Omega}f(u_{\epsilon}-u_0)dx\, ,  \\[0.2cm]
&&\int_{\Omega}A^{ij}_0\frac{\partial u_0}{\partial x_i}\cdot (\frac{\partial u_{\epsilon}}{\partial x_j}-\frac{\partial u_0}{\partial x_j})dx=\int_{\Omega}f(u_{\epsilon}-u_0)dx \, .
\end{eqnarray*}
Subtracting yields
\begin{equation}\label{fundamental-equality}
\int_{\Omega}A^{ij}_{\epsilon}(\frac{\partial u_{\epsilon}}{\partial x_i}-\frac{\partial u_0}{\partial x_i})(\frac{\partial u_{\epsilon}}{\partial x_j}-\frac{\partial u_0}{\partial x_j})dx=\int_{\Omega}(A^{ij}_0-A^{ij}_{\epsilon})\frac{\partial u_0}{\partial x_i}(\frac{\partial u_{\epsilon}}{\partial x_j}-\frac{\partial u_0}{\partial x_j})dx \, ,
\end{equation}
hence, by (\ref{0}) and H\"{o}lder's inequality we obtain
\[
\|\nabla u_{\epsilon}-\nabla u_0\|_{L^2(\Omega)} \leq c\|\nabla u_0\|_{L^q(\Omega)} \|A_{\epsilon}-A_0\|_{L^{2q/(q-2)}(\Omega)}.
\]

{\em Part (2).} Let us denote by $\nabla$ the gradient operator, from $H^1_0(\Omega)$ to $L^2(\Omega)$. By \nl
\cite[Lemma 3]{Ba} we have
\[
u_{\epsilon}-u_0 =\nabla^* A_{\epsilon}^{1/2}(A_{\epsilon}^{1/2}\nabla\nabla^* A_{\epsilon}^{1/2} +1)^{-1}A_{\epsilon}^{-1/2}(A_0-A_{\epsilon})A_0^{-1/2}
(A_0^{1/2}\nabla\nabla^*A_0^{1/2}+1)^{-1}A_0^{1/2}\nabla f,
\]
therefore
\begin{eqnarray*}
\|\nabla u_{\epsilon}-\nabla u_0\|_{L^2(\Omega)} &\leq& \|A_{\epsilon}^{-1/2}\|\cdot\|A_{\epsilon}^{1/2}\nabla \nabla^* A_{\epsilon}^{1/2}(A_{\epsilon}^{1/2}\nabla\nabla^* A_{\epsilon}^{1/2} +1)^{-1}\|  \times \\
&&   \| A_{\epsilon}^{-1/2}(A_0-A_{\epsilon})A_0^{-1/2}
(A_0^{1/2}\nabla\nabla^*A_0^{1/2}+1)^{-1}A_0^{1/2}\nabla \| \cdot \|f\|_{L^2(\Omega)}\\
&\leq & c \| (A_0-A_{\epsilon})A_0^{-1/2}(A_0^{1/2}\nabla\nabla^*A_0^{1/2}+1)^{-1}A_0^{1/2}\nabla \|\, ,
\end{eqnarray*}
where non-indexed norms are operator norms on $L^2(\Omega)$.
By \cite[Lemma 4]{Ba} we conclude that
\begin{equation}
\|\nabla u_{\epsilon}-\nabla u_0\|_{L^2(\Omega)} \leq c\| A_{\epsilon}-A_0\|_{L^r(\Omega)}  \, ,
\mbox{ all }r>\frac{2q}{q-2}(\frac{N}{2}+2\gamma-1),
\label{ch}
\end{equation}
and the result follows.

{\em Part (3).}  Going back to \eqref{fundamental-equality}, applying \eqref{0} to the left-hand side and decomposing the 
integral in the right-hand side by integrating in $K$ and in $\Omega\setminus K$, we get by H\"older inequality for any $q>2$, 
\begin{eqnarray*}
\frac{1}{c}\|\nabla u_\eps-\nabla u_0\|_{ L^2(\Omega)}^2 \leq C\|A_\eps-A_0\|_{L^\infty(\Omega\setminus K)}\|\nabla u_\eps-\nabla u_0\|_{ L^2(\Omega)}+ \\
C\|A_\eps-A_0\|_{L^{2q/(q-2)}(K)}\|\nabla u_0\|_{L^q(K)}\|\nabla u_\eps-\nabla u_0\|_{ L^2(\Omega)}
\end{eqnarray*}
which implies,
\begin{eqnarray*}
\|\nabla u_\eps-\nabla u_0\|_{ L^2(\Omega)} \leq \tilde C(\|A_\eps-A_0\|_{L^\infty(\Omega\setminus K)}+\|A_\eps-A_0\|_{L^{2q/(q-2)}(K)}\|\nabla u_0\|_{L^q(K)})
\end{eqnarray*}

But interior estimates for the limit problem imply, see \cite[Theorem 2]{Me}, that there exists a $q>2$ such that $\|\nabla u_0\|_{L^q(K)}\leq C$. Hence, we obtain the stated result.

{\em Part (4).}  Let us start noticing that in general, since $u_\eps\to u_0$ weakly in $H^1(\Omega)$ and in particular in $L^2(\Omega)$, we
get

\begin{equation}\label{eq-convergence}
\int_{\Omega} A^{ij}_\eps\frac{\partial u_\eps}{\partial x_i}\frac{\partial u_\eps}{\partial x_j}=\int_{\Omega} fu_\eps \to \int_{\Omega}f u_0=\int_{\Omega} A^{ij}_0\frac{\partial u_0}{\partial x_i}\frac{\partial u_0}{\partial x_j}
\end{equation}

Hence we have
\begin{eqnarray*}
\int_{\Omega}A^{ij}_0\frac{\partial u_\eps}{\partial x_i}\frac{\partial u_\eps}{\partial x_j} dx&\leq&
\int_{\Omega}(A_0-A_{\epsilon})_+^{ij}\frac{\partial u_\eps}{\partial x_i}\frac{\partial u_\eps}{\partial x_j}dx +\int_{\Omega}A^{ij}_{\epsilon}\frac{\partial u_\eps}{\partial x_i}\frac{\partial u_\eps}{\partial x_j}dx\\
&\to& \int_{\Omega}A^{ij}_0\frac{\partial u_0}{\partial x_i}\frac{\partial u_0}{\partial x_j} dx \, \qquad \mbox{ as }\epsilon\to 0.
\end{eqnarray*}
where we have used the hypothesis and \eqref{eq-convergence}. Hence
\[
\limsup_{\epsilon\to 0}\int_{\Omega}(A^{ij}_0\frac{\partial u_\eps}{\partial x_i}\frac{\partial u_\eps}{\partial x_j} +|u_{\epsilon}|^2)dx 
\leq \int_{\Omega}(A_0^{ij}\frac{\partial u_0}{\partial x_i}\frac{\partial u_0}{\partial x_j}  + |u_0|^2)dx\, .
\]
Let us now consider the space $H^1(\Omega)$ with the Hilbert norm
\[
\triplenorm{u}^2 =\int_{\Omega}(A^{ij}_0\frac{\partial u}{\partial x_i}\frac{\partial u}{\partial x_j} +|u|^2)dx .
\]
We have thus proved that $\limsup \triplenorm{u_{\epsilon}}\leq\triplenorm{u_0}$. Moreover, the weak convergence in
$H^1_0(\Omega)$ with respect to the standard norm implies the weak convergence with respect to the equivalent norm $\triplenorm{\cdot}$, hence $\liminf \triplenorm{u_{\epsilon}}\geq\triplenorm{u_0}$. We thus conclude that $\triplenorm{u_{\epsilon}}\to \triplenorm{u_0}$, and therefore
$u_{\epsilon}\to u_0$ strongly in $H^1_0(\Omega)$, as required. $\hfill\Box$

{\bf Remark. } If $\Omega$ has a Lipschitz boundary, in the sense that there exists a bi-Lipschitz transformation that maps $\Omega$ onto a domain with $C^1$ boundary, then condition (2) is satisfied. This follows easily from \cite[Theorem 1]{Me}.

{\bf Remark.} It is an interesting problem to prove or disprove that $u_{\epsilon}\to u_0$ in $H^1_0(\Omega)$ without any assumptions other than those stated before Theorem \ref{thm}.

We now present an example that shows that estimate (\ref{est}) is sharp. We fix an angle $\pi<\beta <2\pi$ and denote by $\Omega_{\beta}$ the circular sector of radius one and angle $\beta$,
\[
\Omega_{\beta}=\{ (r,\theta) \; : \; 0<r< 1 \; , \;\; 0<\theta<\beta\}.
\]
Let $u_0$ be the solution of the problem,
\begin{equation}
\label{lim}
\darr{-\Delta u_0=\frac{4\beta^2-\pi^2}{\beta^2}\sin(\frac{\pi\theta}{\beta}),}{\mbox{ in }\Omega_{\beta},}{u_0=0,}{\mbox{ on }\po_{\beta}.}
\end{equation}
(The factor $(4\beta^2-\pi^2)/\beta^2$ is introduced to simplify subsequent calculations.) Simple computations give
\begin{equation}
u_0(x)=(r^{\frac{\pi}{\beta}}-r^2)\sin(\frac{\pi\theta}{\beta}) \; , \quad 0<r<1, \;\;\; 0<\theta<\beta\, .
\label{u0}
\end{equation}
This implies in particular that $\nabla u_0\in L^q(\Omega_{\beta})$ if and only if $q<2\beta/(\beta-\pi)$.
We now fix $\alpha>0$ and for $0<\epsilon<1$ we set
\[
a_{\epsilon}(x)=\darr{\alpha,}{0<|x|<\epsilon,}{1,}{\epsilon<|x|<1.}
\]
We note that $a_{\epsilon}\to 1$ in $L^p(\Omega_{\beta})$, $1\leq p<\infty$, but not in $L^{\infty}(\Omega_{\beta})$. 
Let us denote by $u_{\epsilon}$ the solution of the perturbed boundary-value problem
\begin{equation}
\label{11}
\darr{-\diver(a_{\epsilon}\nabla u_{\epsilon})=\frac{4\beta^2-\pi^2}{\beta^2}\sin(\frac{\pi\theta}{\beta}),}{\mbox{ in }\Omega_{\beta},}{u_{\epsilon}=0,}{\mbox{ on }\po_{\beta}.}
\end{equation}
From part (2) of Theorem \ref{thm} follows that
\[
\|\nabla u_{\epsilon}- \nabla u_0\|_{L^2(\Omega_{\beta})}\leq c\|a_{\epsilon}-a_0\|_{L^p(\Omega_{\beta})}=c\epsilon^{\frac{2}{p}}
 \; , \qquad \mbox{ all }p>\frac{2\beta}{\pi}.
\]
Now, simple computations give that $u_{\epsilon}(x) =v_{\epsilon}(r)\sin(\pi\theta/\beta)$, where
\[
v_{\epsilon}(r)=\darr{\Frac{(\frac{1}{\alpha}-1)(\epsilon^{2}+\epsilon^{2-\frac{2\pi}{\beta}})+2\epsilon^{-\frac{\pi}{\beta}}}
{(1-\alpha)\epsilon^{\frac{\pi}{\beta}}+(1+\alpha)\epsilon^{-\frac{\pi}{\beta}}}r^{\frac{\pi}{\beta}}-\frac{1}{\alpha}r^2,}{0<r<\epsilon,}
%%%
{\Frac{(1-\alpha)\epsilon^2+(1+\alpha)\epsilon^{-\frac{\pi}{\beta}}}{(1-\alpha)\epsilon^{\frac{\pi}{\beta}}+(1+\alpha)\epsilon^{-\frac{\pi}{\beta}}}
r^{\frac{\pi}{\beta}}
+\frac{(1-\alpha)(\epsilon^{\frac{\pi}{\beta}}-\epsilon^2)}{(1-\alpha)\epsilon^{\frac{\pi}{\beta}}+
(1+\alpha)\epsilon^{-\frac{\pi}{\beta}}}r^{-\frac{\pi}{\beta}}-r^2,}{\epsilon<r<1.}
\]
Let us denote $u_0(x)=v_0(r)\sin(\pi\theta/\beta)$ (cf. (\ref{u0})). We then have
\begin{eqnarray*}
\| \nabla u_{\epsilon}-\nabla u_0\|_{L^2(\Omega_{\beta})}^2 &\geq& c\int_0^{\epsilon}(v'_{\epsilon}-v'_0)^2 r\, dr\\
&=&\frac{\pi(1-\alpha)^2}{2\beta(1+\alpha)^2}\epsilon^{\frac{2\pi}{\beta}}(1+o(1)) \; , \qquad \mbox{ as }\epsilon \to 0.
\end{eqnarray*}
We conclude that the index $2q/(q-2)$ in (\ref{est}) cannot be replaced by any smaller index.

\

\section{Stability estimates in $H^1_0$ under domain perturbation}

Let $\Omega$ be a bounded domain in $\R^N$. We fix $M>0$. Let $\phi:\Omega\to \phi(\Omega)=:\tilde\Omega$ be a bi-Lipschitz transformation such that
\begin{equation}
\|D\phi\|_{L^{\infty}(\Omega)}\leq M , \qquad  \|(D\phi)^{-1}\|_{L^{\infty}(\Omega)}\leq M \; .
\label{s1}
\end{equation}
We set $E=\{x\in\Omega \; : \; \phi(x)\neq x\}$ and we note that
\begin{equation}
\| F\circ\phi -F\|_{L^2(\Omega)}\leq \|F\|_{L^q(\Omega)} |E|^{\frac{q-2}{2q}}  \; , \qquad
\| (D\phi)^{\pm 1} -I\|_{L^q(\Omega)} \leq c |E|^{\frac{q-2}{2q}}.
\label{wwww}
\end{equation}
Let $A=\{A_{ij}(x)\}$ be measurable and symmetric on $\Omega\cup\tilde\Omega$ satisfying
\[
\frac{1}{M}|\xi|^2 \leq \sum_{i,j}A_{ij}(x)\xi_i\xi_j \leq M |\xi|^2 , \qquad  x\in \Omega\cup\tilde\Omega \; , \;\;  \; \xi\in\R^N.
\]
Let $L$ (resp. $\tilde L$) denote the operator $-\sum_{i,j}\partial_{x_j}\{A_{ij}\partial_{x_i}\}$ on $L^2(\Omega)$ (resp. $L^2(\tilde\Omega)$) subject to
Dirichlet boundary conditions. We fix $f\in L^2(\Omega\cup\tilde\Omega)$ and we denote by $u$ (resp. $\tilde u$) the solution of
$Lu = f$  (resp. $\tilde L\tilde u=f$). 
We extend $\nabla u$ and $\nabla \tilde u$ to be zero outside their respective domains. We then have the following
%%%%%%%%%%%%%%%%%%%%%%%%%%%%%%%%%%%%%%%%%%%%%%%%%%%%%
\begin{theorem}
Assume that there exists $q>2$ such that $\|f\|_{L^q(\Omega\cup\tilde\Omega)}\leq M$,
$\|\nabla u\|_{L^q(\Omega)}\leq M$ and $\|\nabla \tilde u\|_{L^q(\tilde\Omega)}\leq M$. Then there exists a constant $c$ depending only on $M$ such that
\begin{equation}
\|\nabla \tilde u -\nabla u\|_{L^2(\Omega\cup\tilde\Omega)} \leq c|E|^{\frac{q-2}{2q}} .
\label{eq:conclusion}
\end{equation}
\label{thm:stab_poiss}
\end{theorem}
%%%%%%%%%%%%%%%%%%%%%%%%%%%%%%%%%%%%%%%%%%%%%%%%%%%%%
{\em Proof.} We set $g(x)=|\det D\phi(x)|$, $x\in\Omega$. We note that given
$v\in H^1_0(\tilde\Omega)$ and letting $u=v\circ \phi$ we have
\[
\int_{\tilde\Omega}|v|^2dy =\int_{\Omega}|u|^2 g \, dx\; 
\]
and
\[
\int_{\tilde\Omega}   \sum_{i,j=1}^NA_{ij}\parder{v}{y_i}\parder{\bar v}{y_j} dy =\int_{\Omega}\sum_{i,j=1}^N a_{ij}\parder{u}{x_i}\parder{\bar u}{x_j}g  \, dx\; ,
\]
where $a=(a_{ij})_{i,j=1, \dots , N}$ is defined on $\Omega$ by
\begin{equation}
a_{ij}=\sum_{r,s=1}^N \Big(A_{rs}  \parder{  \phi_{i}^{(-1)}}{y_r}   \parder{\phi_{j}^{(-1)}}{y_s}\Big)\circ \phi = ((D\phi)^{-1} (A\circ\phi)  (D\phi)^{-t})_{ij}\, .
\label{matrixa}
\end{equation}
This leads to the notion of the {\em pull-back}: the pull-back of $\tilde L$ to $\Omega$ is the self-adjoint operator $\tilde H$ on $L^2(\Omega,g\, dx)$
associated with the sesquilinear form $\tilde Q$ with $\dom(\tilde Q) = H^1_0(\Omega)$  and
\[
\tilde Q(u_1,u_2)=\int_{\Omega}\sum_{i,j=1}^N a_{ij}\parder{u_1}{x_i}\parder{\bar u_2}{x_j}g \, dx, \ \qquad u_1,u_2\in 
H^1_0(\Omega).
\]
So formally
$
\tilde H U= -g^{-1}\sum_{i,j}( g a_{ij}U_{x_i})_{x_j}$.
Equivalently, $\tilde H$ can be defined as $\tilde H=C_{\phi}\tilde L C_{\phi}^{-1}$, where
$C_{\phi}:L^2(\tilde\Omega)\to L^2(\Omega, g\, dx)$ denotes composition by $\phi$ (a unitary operator).

It then follows that $\tilde H=g^{-1}\hat{H}$, where $\hat{H}$ is the self-adjoint operator on $L^2(\Omega)$ associated to the form $\tilde Q$ defined above on $H^1_0(\Omega)$. Now, let $\hat{u}\in L^2(\Omega)$ be defined by
$\hat{H}\hat{u} =f$. Using (1) of Theorem \ref{thm} and also (\ref{wwww}) we obtain
\begin{eqnarray*}
\|\nabla\hat{u} -\nabla u\|_{L^2(\Omega)}&\leq &c\| a - A\|_{L^{\frac{2q}{q-2}}(\Omega)} \\
&\leq& c |E|^{\frac{q-2}{2q}}.
\end{eqnarray*}
Morover we have
\begin{equation}
\darr{\tilde H (\tilde u \circ\phi)=f\circ\phi,}{}
{\tilde H\hat{u}=gf.}{}
\label{kath}
\end{equation}
Hence, by (\ref{wwww}),
\begin{equation}
\| \nabla (\tilde u\circ\phi) -\nabla \hat{u}\|_{L^2(\Omega)}
\leq c\|f\circ\phi -gf\|_{L^2(\Omega)} \leq c |E|^{\frac{q-2}{2q}}.
\label{cof}
\end{equation}
Finally by (\ref{wwww}) we also have
\begin{equation}
\| \nabla (\tilde u\circ\phi) -\nabla \tilde u\|_{L^2(\Omega)}\leq c|E|^{\frac{q-2}{2q}}.
\label{ppp}
\end{equation}
Combining (\ref{kath}), (\ref{cof}) and (\ref{ppp}) we obtain that
$\|\nabla \tilde u -\nabla u\|_{L^2(\Omega)} \leq c|E|^{\frac{q-2}{2q}}$.
We also have
\[
|\tilde \Omega\setminus\Omega| \leq c|\phi^{-1}(\tilde\Omega\setminus\Omega)| \leq c|E|,
\]
hence $\|\nabla \tilde u -\nabla u\|_{L^2(\tilde\Omega)} \leq c|E|^{\frac{q-2}{2q}}$ by H\"{o}lder inequality.
This concludes the proof. $\hfill\Box$

Using Theorem \ref{thm:stab_poiss} we can now prove stability estimates in $H^1_0$ for localized perturbations in terms of the Lebesgue measure. We consider the following class of domains (see also \cite{BL1}):

\begin{definition}
Let $V$ be a bounded open cylinder, i.e., there exists a rotation  $R$  such that $R(V)=W\times ]a,b[$, where $W$ is a bounded convex open set in $\R^{N-1}$. Let $M,\rho>0$. We say that a bounded open set $\Omega\subset\R^N$ belongs to $\cC^{0,1}_M(V,R)$ if
$\Omega$ is of class $C^{0,1}$ (i.e., $\Omega$ is locally a subgraph of $C^{0,1}$ functions) and there exists a function $h\in C^{0,1}(\overline{W})$ such that $a+(b-a)/10 \leq h\le b$,   ${\rm Lip}(h)\le M$, and
\begin{equation}
R(\Omega\cap V)=\{ (\bar{x},x_N) \; : \;  \bar{x}\in W \, ,\,  a<x_N<h(\bar{x}) \}.
\label{C11}
\end{equation}
\label{def:C11}
\end{definition}
In the following theorem, again, $\nabla u$ and $\nabla\tilde u$ are extended to be zero outside $\Omega$ and $\tilde\Omega$ respectively.
\begin{theorem}
Let $\Omega$, $\tilde{\Omega}$ be domains of class $\cC^{0,1}_M(V,R)$ with $\Omega\setminus V=\tilde\Omega\setminus V$. Let $f\in L^2(\Omega\cup\tilde\Omega)$ and let $u\in H^1_0(\Omega)$ and $\tilde u\in H^1_0(\tilde\Omega)$ solve $Lu=f$ and $\tilde L \tilde u=f$ respectively. Assume that there exists $q>2$ such that $\|f\|_{L^q(\Omega\cup\tilde\Omega)}\leq M$,
$\|\nabla u\|_{L^q(\Omega)}\leq M$ and $\|\nabla \tilde u\|_{L^q(\tilde\Omega)}\leq M$. Then there exists a constant $c$ depending only on $M$ such that
\begin{equation}
\|\nabla \tilde u -\nabla u\|_{L^2(\Omega\cup\tilde\Omega)} \leq c|\tilde\Omega\triangle\Omega|^{\frac{q-2}{2q}} \; .
\label{fuga}
\end{equation}
\label{thm:niki}
\end{theorem}
{\em Proof.} It has been proved in \cite{BBL}, following earlier work of Burenkov and Lamberti \cite{BL1}, that if $\Omega$ and $\tilde{\Omega}$ are
in $\cC^{0,1}_M(V,R)$ then there exists a bi-Lipschitz diffeomorphism $\phi:\Omega\to\tilde\Omega$ whose Lipschitz constants are estimated in terms of $M$ only and such that $|E| \leq c |\tilde\Omega\triangle\Omega|$, where $c$ also depends only on $M$. Hence (\ref{fuga}) follows from Theorem \ref{thm:stab_poiss}. $\hfill\Box$

We now present an example that shows that estimate (\ref{est}) is sharp. It is a variation of the example in Section \ref{sec2}.
We fix $\pi<\beta <2\pi$ and $0\leq\epsilon<1$ and denote by $\Omega_{\epsilon, \beta}$ the planar domain
\[
\Omega_{\epsilon, \beta}=\{ (r,\theta) \; : \; \epsilon <r< 1 \; , \;\; 0<\theta<\beta\}.
\]
We denote by $u_{\epsilon}$ be the solution of the problem
\begin{equation}
\darr{-\Delta u_{\epsilon}=\frac{4\beta^2-\pi^2}{\beta^2}\sin(\frac{\pi\theta}{\beta}),}{\mbox{ in }\Omega_{\epsilon, \beta},}{u_{\epsilon}=0,}{\mbox{ on }\po_{\epsilon, \beta}.}
\end{equation}
The function $u_0$ has been computed in (\ref{u0}). For $\epsilon>0$ a direct computation gives that $u_{\epsilon}(x)=v_{\epsilon}(r)\sin(\pi\theta/\beta)$, where
\[
v_{\epsilon}(r)=\frac{\epsilon^{-\pi/\beta}-\epsilon^2}{\epsilon^{-\pi/\beta}-\epsilon^{\pi/\beta}}r^{\frac{\pi}{\beta}} +
\frac{\epsilon^2 -\epsilon^{\pi/\beta}}{\epsilon^{-\pi/\beta}-\epsilon^{\pi/\beta}}r^{-\frac{\pi}{\beta}}  -r^2 \; , \qquad \epsilon<r<1.
\]
We then easily obtain the asymptotic formula
\begin{equation}
\label{actual}
\|\nabla u_{\epsilon}-\nabla u_0\|_{L^2(\Omega_0)}  =A\epsilon^{\pi/\beta}(1+o(1)),
\end{equation}
for some $A>0$.

Let us now see what our theorem gives. We note that the domains are not of some class $\cC^{0,1}_M(V,R)$ uniformly in $\epsilon$, so we cannot directly apply Theorem \ref{thm:niki}. Instead, let $s_{\epsilon}:(0,1)\to (\epsilon,1)$,
\[
s_{\epsilon}(r)= \darr{\frac{r}{2}+\epsilon \, ,}{0<r<2\epsilon,}{r,}{2\epsilon<r<1.}
\]
We define the bi-Lipschitz map $\phi_{\epsilon}:\Omega_0\to \Omega_{\epsilon}$, given in polar coordinates by
\[
\phi_{\epsilon}(r,\theta)=(s_{\epsilon}(r),\theta) \; , \qquad 0<r<1 \; , \qquad 0<\theta<\beta.
\]
Using the explicit computation of $u_0$ and $u_{\epsilon}$ we easily see that the assumptions of Theorem \ref{thm:stab_poiss} are satisfied for any (fixed) $q<2\beta/(\beta-\pi)$. Hence Theorem \ref{thm:stab_poiss} gives
\[
\|\nabla u_{\epsilon} -\nabla u_0\|_{L^2(\Omega_0\cup\Omega_{\epsilon})} \leq c_q|E_{\epsilon}|^{\frac{q-2}{2q}}=c_q\epsilon^{\frac{q-2}{q}}
=c_{\delta}\epsilon^{\frac{\pi}{\beta}-\delta},
\]
for any $\delta>0$. Because of (\ref{actual}), this shows that the exponent of $|E|$ in (\ref{eq:conclusion}) cannot be replaced by a smaller one.

%%%%%%%%%%%%%%%%%%%%%%%%%%%%%%%%%%%%%%%%%%%%%%%%%%%%%%%%%%%%%%%%%%%%%%%%%%%%%%%%%%%%%

{\small
\noindent Jos\'{e} M. Arrieta\\
Departamento de Matem\'{a}tica Aplicada\\
Universidad Complutense de Madrid\\
28040 Madrid\\
Spain\\
e-mail: arrieta@mat.ucm.es

\noindent Gerassimos Barbatis\\
Department of Mathematics\\
University of Athens\\
15784 Athens\\
Greece\\
e-mail: gbarbatis@math.uoa.gr

\noindent 
}

\end{document}